\documentclass[12pt,twoside]{article}
\usepackage{amssymb}

\font\teneufm=eufm10 scaled \magstep1
\font\seveneufm=eufm7 scaled \magstep1
\font\fiveeufm=eufm5  scaled \magstep1
\newfam\eufmfam
\textfont\eufmfam=\teneufm
\scriptfont\eufmfam=\seveneufm
\scriptscriptfont\eufmfam=\fiveeufm
\def\frak#1{{\fam\eufmfam\relax#1}}

\newfam\msbfam
\font\tenmsb=msbm10 scaled \magstep1  \textfont\msbfam=\tenmsb
\font\sevenmsb=msbm7 scaled \magstep1 \scriptfont\msbfam=\sevenmsb
\font\fivemsb=msbm5 scaled \magstep1  \scriptscriptfont\msbfam=\fivemsb
\def\Bbb{\fam\msbfam \tenmsb}

\def\RR{{\Bbb R}}
\def\CC{{\Bbb C}}
\def\QQ{{\Bbb Q}}
\def\NN{{\Bbb N}}
\def\ZZ{{\Bbb Z}}

\def\PP{{\Bbb P}}

\def\ra{\rightarrow}

 \def\HollowBoxx #1#2#3{{\dimen0=#1 \advance\dimen0 by -#2
       \dimen1=#1 \advance\dimen1 by #3
        \vrule height 0pt depth #3 width #2
       \hskip -#3
       \vrule height #1 depth #3 width #3}}
 \def\LeftContraction{\mathord{\kern1.45pt \HollowBoxx{6pt}{3.5pt}{.4pt}}\,}

 \def\HollowBox #1#2#3{{\dimen0=#1 \advance\dimen0 by -#3
       \dimen1=#1 \advance\dimen1 by #3
        \vrule height #1 depth #3 width #3
        \vrule height 0pt depth #3 width #2
        \hskip -#3}}
 \def\RightContraction{\mathord{\, \HollowBox{6pt}{3.1pt}{.4pt}} \kern1.6pt}

\newtheorem{theorem}{THEOREM}[section]

\begin{document}

\begin{center}
{\Large \bf Analogues of Rossi's Map and
\medskip\\
E. Cartan's Classification of
\medskip\\
Homogeneous Strongly Pseudoconvex
\medskip\\ 
3-Dimensional Hypersurfaces}\footnote{{\bf Mathematics Subject Classification:} 53C30, 32V05.}\footnote{{\bf
Keywords and Phrases:} homogeneous 3-dimensional $CR$-manifolds, Rossi's example.}
\medskip\\
\normalsize A. V. Isaev
\end{center}

\begin{quotation} \small \sl We introduce analogues of a map due to Rossi and show how they can be used to explicitly determine all covers of certain homogeneous strongly pseudoconvex 3-dimensional hypersurfaces that appear in the classification obtained by E. Cartan in 1932. 
\end{quotation}

\thispagestyle{empty}

\pagestyle{myheadings}
\markboth{A. V. Isaev}{Analogues of Rossi's Map}

\setcounter{section}{-1}

\section{Introduction}
\setcounter{equation}{0}

In 1932 E. Cartan \cite{C} classified all connected 3-dimensional homogeneous strongly pseudoconvex $CR$-manifolds (we reproduce his classification in detail in Section \ref{classification} below). Most manifolds in the classification are given by explicit equations as hypersurfaces in either $\CC^2$ or $\CC\PP^2$. The exceptions (apart from lens spaces) are all possible covers of each of the following hypersurfaces:
$$
\begin{array}{lll}
\chi & :=&\left\{(z,w)\in\CC^2:x^2+u^2=1\right\},\\
\mu_{\alpha} & := &\left\{(z:w:\zeta)\in\CC\PP^2: |z|^2+|w|^2+|\zeta|^2=\alpha |z^2+w^2+\zeta^2|\right\},\,\alpha>1,\\
\nu_{\alpha} &:=& \left\{(z,w)\in\CC^2: |z|^2+|w|^2-1=\alpha|z^2+w^2-1|\right\}\setminus\\
&&\left\{(x,u)\in\RR^2:x^2+u^2=1\right\},\,-1<\alpha<1,\\
\eta_{\alpha} & := &\Bigl\{(z,w)\in\CC^2: 1+|z|^2-|w|^2=\alpha|1+z^2-w^2|,\\
&&\hbox{Im}(z(1+\overline{w}))>0\Bigr\},\,\alpha>1.
\end{array}
$$
In the above formulas and everywhere below we set $z=x+iy$, $w=u+iv$. 

Sometimes it is desirable to know all the covers explicitly. For example, such a need arises when one attempts to describe all Kobayashi-hyperbolic 2-dimensional complex manifolds with holomorphic automorphism group of dimension 3, since almost every hypersurface from Cartan's classification can be realized as an orbit of the automorphism group action on such a manifold \cite{I2}  (see also \cite{I1} for motivations and related results).

All covers of the hypersurface $\chi$ are easy to find (see Section \ref{chicovers}). Next, in Section \ref{mucovers} we determine the covers of the hypersurfaces $\mu_{\alpha}$. All of them are compact, and, since compact Levi non-degenerate homogeneous $CR$-manifolds have been extensively studied (see, e.g., \cite{AHR} and references therein), all these covers have been discovered on many occasions in various guises (see, e.g. \cite{ENS}). In Section \ref{mucovers} we show how the covers of $\mu_{\alpha}$ can be found by using a map introduced by Rossi in \cite{R}. Originally, Rossi defined it as a map from $\CC\PP^2\setminus\{0\}$ into $\CC\PP^3$ and used it to construct strongly pseudoconvex 3-dimensional manifolds that do not bound any pseudoconvex analytic space. For our purposes, however, we will only be interested in the restriction of the map to $\CC^2\setminus\{0\}$ and call this restriction {\it Rossi's map} (see formula (\ref{rossioriginal})).

Our main results are contained in Section \ref{coversnueta}, where we determine all covers of the hypersurfaces $\nu_{\alpha}$ and $\eta_{\alpha}$ that, in contrast with those of $\mu_{\alpha}$, appear to have never been found explicitly. Interestingly, it turns out that this can be done by introducing a map (in fact, a sequence of maps) analogous to Rossi's map (see (\ref{analoguerossi}), (\ref{rossin})). While Rossi's map is associated with the action of $SU_2$ on $\CC^2\setminus\{0\}$, our analogues are associated with the action of $SU_{1,1}$ on $\CC^2\setminus\{(z,w)\in\CC^2: |z|^2-|w|^2=0\}$ (which is a manifestation of the non-compactness of $\nu_{\alpha}$ and $\eta_{\alpha}$). While for $\mu_{\alpha}$ only 2- and 4-sheeted covers can occur, there is an $n$-sheeted cover for every $n\ge 2$ as well as an infinite-sheeted cover, for each of $\nu_{\alpha}$, $\eta_{\alpha}$. In fact, as in the cases of $\chi$ and $\mu_{\alpha}$, we explicitly find the covers of certain domains in $\CC^2$, namely $D^{\nu}:=\cup_{-1<\alpha<1}\nu_{\alpha}$ and $D^{\eta}:=\cup_{\alpha>1}\eta_{\alpha}$. We equip the covers of these domains with complex structures obtained by pull-backs under explicit covering maps, the covers of $\nu_{\alpha}$, $\eta_{\alpha}$ with $CR$-structures induced by these complex structures, and find the groups of $CR$-automorphisms of the induced $CR$-structures explicitly. Hence the emphasis of the present paper is on the explicit determination of the covers of $\nu_{\alpha}$ and $\eta_{\alpha}$ in the differential-geometric sense, and on group actions.

The main motivation for the present work was the above-mentioned problem of classifying 2-dimensional Kobayashi-hyperbolic manifolds with\linebreak 3-dimensional automorphism groups (we call such manifolds (2,3)-{\it manifolds} for brevity). A complete explicit classification of (2,3)-manifolds was obtained in \cite{I2}, and the particular realizations of the covers of  $\nu_{\alpha}$ and $\eta_{\alpha}$ (as well as those of $\chi$ and $\mu_{\alpha}$) constructed in the present paper turned out to be extremely useful for this purpose. In our construction, the covers $\nu_{\alpha}^{(N)}$ and $\eta_{\alpha}^{(N)}$ of $\nu_{\alpha}$ and $\eta_{\alpha}$ with the same number of sheets $N$ (where $N$ can be infinite) are glued together into complex manifolds $M_N^{\nu}$ and $M_N^{\eta}$, respectively. On these manifolds certain connected 3-dimensional Lie groups $G_N^{\nu\,c}$ and $G_N^{\eta\,c}$ act by holomorphic transformations (these groups are in fact the connected identity components of the automorphism groups of $M_N^{\nu}$ and $M_N^{\eta}$), and the orbits of these actions are exactly the hypersurfaces $\nu_{\alpha}^{(N)}$ and $\eta_{\alpha}^{(N)}$, respectively (the groups  $G_N^{\nu\,c}$ and $G_N^{\eta\,c}$ also coincide with the connected identity components of the groups of $CR$-automorphisms of $\nu_{\alpha}^{(N)}$ and $\eta_{\alpha}^{(N)}$) -- see Theorem \ref{main}. The actions of $G_N^{\nu\,c}$ and $G_N^{\eta\,c}$ on $M_N^{\nu}$ and $M_N^{\eta}$ are proper without singular orbits, hence the orbit spaces $M_N^{\nu}/G_N^{\nu\,c}$ and $M_N^{\eta}/G_N^{\eta\,c}$ are homeomorphic to $\RR$ (see \cite{AA}). The manifolds $M_N^{\nu}$ and $M_N^{\eta}$ cover the domains $D^{\nu}$ and $D^{\eta}$ by means of $N$-to-1 explicit covering maps ${\frak P}_N^{\nu}$ and ${\frak P}_N^{\eta}$ that are expressed in terms of our analogues of Rossi's map. The restrictions of ${\frak P}_N^{\nu}$ to $\nu_{\alpha}^{(N)}$ and ${\frak P}_N^{\eta}$ to $\eta_{\alpha}^{(N)}$ are $N$-to-1 covering maps onto $\nu_{\alpha}$ and $\eta_{\alpha}$, respectively. Thus, our construction immediately gives a large number of examples of 
(2,3)-manifolds whose automorphism groups act with codimension 1 orbits: indeed, any domain in $M_N^{\nu}$ bounded by two $G_N^{\nu\,c}$-orbits and any domain in $M_N^{\eta}$ bounded by two $G_N^{\eta\,c}$-orbits are such manifolds (this phenomenon takes place for other types of hypersurfaces from Cartan's classification as well).

Another feature of our construction utilized in \cite{I2} concerns the complex structures of $M_N^{\nu}$ and $M_N^{\eta}$. These structures are given by pull-backs under the maps ${\frak P}_N^{\nu}$ and ${\frak P}_N^{\eta}$ and, since the formulas for ${\frak P}_N^{\nu}$ and ${\frak P}_N^{\eta}$ are complicated, the resulting complex structures are not very explicit. Nevertheless, some properties of these structures can be derived from our construction. In \cite{I2}, in order to classify 
(2,3)-manifolds whose automorphism groups act with codimension 2 orbits, we study complex curves in $M_N^{\nu}$ and $M_N^{\eta}$ invariant under the actions of maximal compact subgroups of $G_N^{\nu\,c}$ and $G_N^{\eta\,c}$, respectively, for $N<\infty$. Namely, let $K\subset G_N^{\nu\,c}$ be a maximal compact subgroup (every such subgroup is isomorphic to the circle) and $C$ a connected $K$-invariant non-singular complex curve in $M_N^{\nu}$. Next, let $G^{\nu}$ be the connected identity component of the automorphism group of $D^{\nu}$ (in fact, $G^{\nu}$ is isomorphic to the connected identity component $SO_{2,1}(\RR)^c$ of $SO_{2,1}(\RR)$). The group $G_N^{\nu\,c}$ covers the group $G^{\nu}$ by means of an explicit $N$-to-1 covering homomorphism ${\cal P}_N^{\eta}$, the group $K':={\cal P}_N^{\nu}(K)$ is a maximal compact subgroup of $G^{\nu}$, and the restriction of ${\cal P}_N^{\nu}$ to $K$ is an $N$-to-1 covering homomorphism onto $K'$. Then it turns out that $C':={\frak P}_N^{\nu}(C)$ is a $K'$-invariant non-singular complex curve in $D^{\nu}$, and the restriction of ${\frak P}_N^{\nu}$ to $C$ is an $N$-to-1 covering map onto $C'$. Moreover, every connected non-singular $K'$-invariant complex curve in $D^{\nu}$ is obtained in this way. A similar statement holds for complex curves in $M_N^{\eta}$ invariant under the actions of maximal compact subgroups of $G_N^{\eta\,c}$. In \cite{I2} we often relied on this property when studying (2,3)-manifolds with codimension 2 orbits.       

Before proceeding, we would like to thank one of the anonymous referees for numerous suggestions that help improve the paper. These suggestions included, in particular, a group-theoretic interpretation of Rossi's map that we reproduce in Section \ref{mucovers}. Similar constructions lead to group-theoretic interpretations of our analogues of Rossi's map (see Section \ref{coversnueta}). 

\section{E. Cartan's Classification}\label{classification}
\setcounter{equation}{0}

In this section we reproduce E. Cartan's classification of connected\linebreak 3-dimensional homogeneous strongly pseudoconvex $CR$-manifolds from \cite{C}. He shows that every such hypersurface is $CR$-equivalent to one of the manifolds on the list below. As does Cartan, we group the model manifolds into spherical and non-spherical ones.
\vspace{0cm}\\

\noindent {\bf The Spherical Case}
$$
\begin{array}{ll}
\hbox{(i)} & S^3,\\
\hbox{(ii)}& {\cal L}_m:=S^3/\ZZ_m,\, m\in\NN, n\ge 2,\,\,\hbox{(lens spaces),}\\
\hbox{(iii)} & \sigma:=\left\{(z,w)\in\CC^2:u=|z|^2\right\},\\
\hbox{(iv)} & \sigma_{+}:=\left\{(z,w)\in\CC^2:u=|z|^2,\,x>0\right\},\\
\hbox{(v)} & \varepsilon_{\alpha}:=\left\{(z,w)\in\CC^2: |z|^2+|w|^{\alpha}=1,\, w\ne 0\right\}, \, \alpha>0,\\
\hbox{(vi)} & \omega:=\left\{(z,w)\in\CC^2:|z|^2+e^u=1\right\},\\
\hbox{(vii)} & \delta:=\left\{(z,w)\in\CC^2: |w|=\exp\left(|z|^2\right)\right\},\\
\hbox{(viii)} & \nu_0=S^3\setminus\RR^2,\\
\hbox{(ix)} & \hbox{any cover of $\nu_0$.}
\end{array}
$$

The groups of $CR$-automorphisms of the above hypersurfaces (except in case (ix)) are as follows:
\smallskip\\

$\underline{\hbox{Aut}_{CR}(S^3)\simeq SU_{2,1}/(\hbox{center})}:$
\begin{equation}
\left(
\begin{array}{c}
z\\
w
\end{array}
\right) \mapsto \displaystyle\frac{\left(\begin{array}{cc}
a_{11} & a_{12}\\
a_{21} & a_{22}
\end{array}
\right)\left(
\begin{array}{c}
z\\
w
\end{array}\right)+\left(
\begin{array}{c}
b_1\\
b_2
\end{array}
\right)}{c_1z+c_2w+d},\label{group}
\end{equation}
with $Q\in SU_{2,1}$, where
\begin{equation}
Q:=\left(\begin{array}{ccc}
a_{11} & a_{12} & b_1\\
a_{21} & a_{22} & b_2\\
c_1& c_2 &d
\end{array}
\right).\label{matrix}
\end{equation}
\smallskip\\

$\underline{\hbox{Aut}_{CR}({\cal L}_m)\simeq U_2/\ZZ_m}:$
$$
\left[\left(
\begin{array}{c}
z\\
w
\end{array}
\right)\right] \mapsto
\left[U
\left(
\begin{array}{c}
z\\
w
\end{array}
\right)\right],
$$
where $U\in U_2$ and $[(z,w)]\in{\cal L}_m$ denotes the equivalence class of $(z,w)\in S^3$ under the action of $\ZZ_m$ embedded in $U_2$ as a subgroup of scalar matrices.
\smallskip\\

$\underline{\hbox{Aut}_{CR}(\sigma)\simeq CU_1\ltimes H}:$
$$
\begin{array}{lll}
z & \mapsto & \lambda e^{i\varphi}z+a,\\
w & \mapsto & \lambda^2w+2\lambda e^{i\varphi}\overline{a}z+|a|^2+i\gamma,
\end{array}
$$
where $\lambda\in\RR^*$, $\varphi,\gamma\in\RR$, $a\in\CC$, $CU_1$ denotes the conformal unitary group given by the conditions $a=0$, $\gamma=0$, and $H$ denotes the Heisenberg group given by the conditions $\lambda=1$, $\varphi=0$.
\smallskip\\

$\underline{\hbox{Aut}_{CR}(\sigma_{+})\simeq \RR\ltimes\RR^2}:$
$$
\begin{array}{lll}
z & \mapsto & \lambda z+i\beta,\\
w & \mapsto & \lambda^2w-2i\lambda \gamma z+\gamma^2+i\gamma,
\end{array}
$$
where $\lambda>0$, $\beta,\gamma\in\RR$.
\smallskip\\

$\underline{\hbox{Aut}_{CR}(\varepsilon_{\alpha})}:$
$$
\begin{array}{lll}
z&\mapsto&\displaystyle e^{i\varphi}\frac{z-a}{1-\overline{a}z},\\
\vspace{0mm}&&\\
w &\mapsto&\displaystyle
e^{i\psi}\frac{(1-|a|^2)^{1/\alpha}}{(1-\overline{a}z)^{2/\alpha}}w,
\end{array}
$$
where $\varphi,\psi\in\RR$, $a\in\CC$, $|a|<1$. We have $\hbox{Aut}_{CR}(\varepsilon_{\alpha})\simeq\widetilde{SO}_{2,1}(\RR)^c\times_{\hbox{\tiny loc}} U_1$, if $\alpha\not\in\QQ$, where $\widetilde{SO}_{2,1}(\RR)^c$ is the universal cover of the connected identity component $SO_{2,1}(\RR)^c$ of the group $SO_{2,1}(\RR)$, and $\times_{\hbox{\tiny loc}}$ denotes locally direct product; $\hbox{Aut}_{CR}(\varepsilon_{\alpha})\simeq SO_{2,1}(\RR)^{c\,(n)}\times_{\hbox{\tiny loc}} U_1$, if $\alpha=n/k$, with $n,k\in\ZZ$, $n>0$, $k\ge 0$, $(n,k)=1$, where $SO_{2,1}(\RR)^{c\,(n)}$ is the $n$-sheeted cover of $SO_{2,1}(\RR)^c$.
\smallskip\\    

$\underline{\hbox{Aut}_{CR}(\omega)\simeq\widetilde{SO}_{2,1}(\RR)^c\times_{\hbox{\tiny loc}}\RR}:$
$$
\begin{array}{lll}
z&\mapsto&\displaystyle e^{i\varphi}\frac{z-a}{1-\overline{a}z},\\
\vspace{0mm}&&\\
w &\mapsto&\displaystyle w+2\ln\frac{\sqrt{1-|a|^2}}{1-\overline{a}z}+i\gamma,
\end{array}
$$
where $\varphi,\gamma\in\RR$, $a\in\CC$, $|a|<1$.
\smallskip\\

$\underline{\hbox{Aut}_{CR}(\delta)\simeq U_1\times(U_1\ltimes\RR^2)}:$
$$
\begin{array}{lll}
z & \mapsto & e^{i\varphi}z+a,\\
w & \mapsto &e^{i\psi}\exp\Bigl(2e^{i\varphi}\overline{a}z+|a|^2\Bigr)w,
\end{array}
$$
where $\varphi,\psi\in\RR$, $a\in\CC$.
\smallskip\\

$\underline{\hbox{Aut}_{CR}(\nu_0)\simeq SO_{2,1}(\RR)}:$
This group consists of all maps of the form (\ref{group}) with $Q\in SO_{2,1}(\RR)$, where $Q$ is defined in (\ref{matrix}).
\smallskip\\

\noindent {\bf The Non-Spherical Case}
$$
\begin{array}{ll}
\hbox{(i)} & \tau_{\alpha}:=\left\{(z,w)\in\CC^2:u=x^{\alpha},\,x>0\right\},\\
&\alpha\in(-\infty,-1]\cup(1,2)\cup(2,\infty),\\
\hbox{(ii)} & \xi:=\left\{(z,w)\in\CC^2:u=x\cdot\ln x,\,x>0\right\},\\
\hbox{(iii)} & \chi=\left\{(z,w)\in\CC^2:x^2+u^2=1\right\},\\
\hbox{(iv)} & \hbox{any cover of $\chi$},\\ 
\hbox{(v)} & \rho_{\alpha}:=\left\{(z,w)\in\CC^2:r= e^{\alpha\varphi}\right\},\,\alpha>0,\\
&\hbox{where $(r,\varphi)$ denote the polar coordinates in the $(x,u)$-plane}\\
&\hbox{with $\varphi$ varying from $-\infty$ to $\infty$},\\
\hbox{(vi)} & \mu_{\alpha}=\left\{(z:w:\zeta)\in\CC\PP^2: |z|^2+|w|^2+|\zeta|^2=\alpha |z^2+w^2+\zeta^2|\right\},\,\alpha>1,\\
\hbox{(vii)} & \hbox{any cover of $\mu_{\alpha}$ with $\alpha>1$},\\
\hbox{(viii)} & \nu_{\alpha}=\left\{(z,w)\in\CC^2: |z|^2+|w|^2-1=\alpha|z^2+w^2-1|\right\}\setminus\\
&\left\{(x,u)\in\RR^2:x^2+u^2=1\right\},\,-1<\alpha<1,\,\alpha\ne 0,\\
\hbox{(ix)} & \hbox{any cover of $\nu_{\alpha}$ with $-1<\alpha<1$, $\alpha\ne 0$},\\
\hbox{(x)} & \eta_{\alpha}=\Bigl\{(z,w)\in\CC^2: 1+|z|^2-|w|^2=\alpha|1+z^2-w^2|,\\
&\hbox{Im}(z(1+\overline{w}))>0\Bigr\},\,\alpha>1,\\
\hbox{(xi)} & \hbox{any cover of $\eta_{\alpha}$ with $\alpha>1$}. \end{array}
$$

Below we list the groups of $CR$-automorphisms of the above hypersurfaces excluding cases (iv), (vii), (ix), (xi). Note that it follows from the proof of Lemma 3.3 of \cite{IKru} that the automorphism group of a connected non-spherical homogeneous hypersurface in a 2-dimensional complex manifold has at most two connected components. 
\smallskip\\ 

$\underline{\hbox{Aut}_{CR}(\tau_{\alpha})\simeq \RR\ltimes\RR^2}:$
$$
\begin{array}{lll}
z & \mapsto & \lambda z+i\beta,\\
w & \mapsto & \lambda^{\alpha}w+i\gamma,
\end{array}
$$
where $\lambda>0$, $\beta,\gamma\in\RR$.
\smallskip\\

$\underline{\hbox{Aut}_{CR}(\xi)\simeq\RR\ltimes\RR^2}:$
$$
\begin{array}{lll}
z & \mapsto & \lambda z+i\beta,\\
w & \mapsto & (\lambda\ln\lambda)z+\lambda w+i\gamma,
\end{array}
$$
where $\lambda>0$, $\beta,\gamma\in\RR$.
\smallskip\\

$\underline{\hbox{Aut}_{CR}(\chi)\simeq O_2\ltimes\RR^2}:$
$$
\left(
\begin{array}{c}
z\\
w
\end{array}
\right)\mapsto
A\left(
\begin{array}{c}
z\\
w
\end{array}
\right)+
i\left(
\begin{array}{c}
\beta\\
\gamma
\end{array}
\right),
$$
where $A\in O_2(\RR)$, $\beta,\gamma\in\RR$.
\smallskip\\

$\underline{\hbox{Aut}_{CR}(\rho_{\alpha})\simeq \RR\ltimes\RR^2}:$ 
$$
\left(
\begin{array}{l}
z\\
w
\end{array}
\right)\mapsto e^{\alpha\psi}
\left(
\begin{array}{rr}
\cos\psi & \sin\psi\\
-\sin\psi & \cos\psi
\end{array}
\right)
\left(
\begin{array}{l}
z\\
w
\end{array}
\right)+i
\left(
\begin{array}{l}
\beta\\
\gamma
\end{array}
\right),
$$
where $\psi,\beta,\gamma\in\RR$.
\smallskip\\

$\underline{\hbox{Aut}_{CR}(\mu_{\alpha})\simeq SO_3(\RR)}:$
$$
\left(
\begin{array}{c}
z\\
w\\
\zeta
\end{array}
\right)
\mapsto A
\left(
\begin{array}{c}
z\\
w\\
\zeta
\end{array}
\right),
$$
where $A\in SO_3(\RR)$.
\smallskip\\

$\underline{\hbox{Aut}_{CR}(\nu_{\alpha})\simeq SO_{2,1}(\RR)}:$ This group consists of all maps of the form (\ref{group}) with $Q\in SO_{2,1}(\RR)$, where $Q$ is defined in (\ref{matrix}).
\smallskip\\

$\underline{\hbox{Aut}_{CR}(\eta_{\alpha})\simeq SO_{2,1}(\RR)^c}:$
\begin{equation}
\left(
\begin{array}{c}
z\\
w
\end{array}
\right) \mapsto \displaystyle\frac{\left(\begin{array}{cc}
a_{22} & b_2\\
c_2 & d
\end{array}
\right)\left(
\begin{array}{c}
z\\
w
\end{array}\right)+\left(
\begin{array}{c}
a_{21}\\
c_1
\end{array}
\right)}{a_{12}z+b_1w+a_{11}},\label{groupetaalpha}
\end{equation}
with $Q\in SO_{2,1}(\RR)^c$, where $Q$ is defined in (\ref{matrix}).
\smallskip\\

Thus, to obtain an explicit classification from the above lists, one needs to determine all possible covers of $\chi$, $\mu_{\alpha}$, $\nu_{\alpha}$ (including the spherical hypersurface $\nu_0$), and $\eta_{\alpha}$. 

Let $M$ be an arbitrary manifold, $\tilde M$ its universal cover, $\Pi:\tilde M\ra M$ a covering map and $\Gamma_{\Pi}$ the corresponding group of covering transformations of $\tilde M$. Then an arbitrary manifold that covers $M$ is obtained from $\tilde M$ by factoring it by the action of a subgroup of $\Gamma_{\Pi}$. Hence, in order to find all covers of each of the hypersurfaces $\chi$, $\mu_{\alpha}$, $\nu_{\alpha}$, $\eta_{\alpha}$ we need to determine their universal covers, the corresponding groups of covering transformations and all their subgroups.

\section{The Covers of $\chi$}\label{chicovers}
\setcounter{equation}{0}

\noindent Let $\Phi^{\chi}:\CC^2\ra\CC^2\setminus\{x=0,\,u=0\}$ be the following map:
$$
\begin{array}{lll}
z & \mapsto & e^x\cos y+iu,\\
w & \mapsto & e^x\sin y+iv.
\end{array}
$$
Clearly, $\Phi^{\chi}$ is an infinitely-sheeted covering map. Introduce on the domain of $\Phi^{\chi}$ a complex structure so that $\Phi^{\chi}$ becomes holomorphic (the pull-back complex structure under $\Phi^{\chi}$), and denote the resulting manifold by $M^{\Phi^{\chi}}$. Then $\tilde\chi$ coincides with the hypersurface
$$
\chi^{(\infty)}:=\left\{(z,w)\in M^{\Phi^{\chi}}:x=0\right\},
$$
equipped with the $CR$-structure induced by the complex structure of $M^{\Phi^{\chi}}$. 

Clearly, $\Gamma_{\Phi^{\chi}}$ consists of all transformations of the form
$$
\begin{array}{lll}
z & \mapsto & z+2\pi i k,\quad k\in\ZZ,\\
w & \mapsto & w.
\end{array}
$$
Let $\Gamma\subset\Gamma_{\Phi^{\chi}}$ be a subgroup. Then there exists an integer $n\ge 0$ such that every element of $\Gamma$ has the form
$$
\begin{array}{lll}
z & \mapsto & z+2\pi i n k,\quad k\in\ZZ,\\
w & \mapsto & w.
\end{array}
$$
Suppose that $n\ge 1$ and consider the map $\Phi^{\chi}_n$ from  $\CC^2\setminus\{x=0,u=0\}$ onto itself defined as follows:
$$
\begin{array}{lll}
z & \mapsto & \hbox{Re}\left(x+iu\right)^n+iy,\\
w & \mapsto & \hbox{Im}\left(x+iu\right)^n+iv.
\end{array}
$$
Denote by $M^{\Phi^{\chi}_n}$ the domain of $\Phi^{\chi}_n$ with the pull-back complex structure under $\Phi^{\chi}_n$. Then the hypersurface 
$$
\chi^{(n)}:=\left\{(z,w)\in M^{\Phi^{\chi}_n}:x^2+u^2=1\right\},
$$
equipped with the $CR$-structure induced by the complex structure of $M^{\Phi^{\chi}_n}$ is an $n$-sheeted cover of $\chi$ corresponding to $\Gamma$ with covering map $\chi^{(n)}\ra\chi$ coinciding with $\Phi^{\chi}_n:M^{\Phi^{\chi}_n}:\ra\CC^2\setminus\{x=0,\,u=0\}$ and factorization map $\chi^{(\infty)}\ra\chi^{(n)}$ given by
$$
\begin{array}{lll}
z & \mapsto & e^{x/n}\cos\left(y/n\right)+iu,\\
w & \mapsto & e^{x/n}\sin\left(y/n\right)+iv.
\end{array}
$$ 

Thus, every cover of $\chi$ is $CR$-equivalent to either $\chi^{(
\infty)}$ or $\chi^{(n)}$ for some $n\in\NN$. The groups of $CR$-automorphisms of $\chi^{(\infty)}$ and $\chi^{(n)}$ are given below.
\smallskip\\

$\underline{\hbox{Aut}_{CR}\left(\chi^{(\infty)}\right)\simeq \RR^3\rtimes\ZZ_2}:$ This group is generated by the maps
$$
\begin{array}{ccc}
z & \mapsto & z+i\beta,\\
w & \mapsto & w+a,
\end{array}
$$
where $\beta\in\RR$, $a\in\CC$, that form the connected identity component $\hbox{Aut}_{CR}(\chi^{(\infty)})^c$ of $\hbox{Aut}_{CR}\left(\chi^{(\infty)}\right)$, and the map
$$
\begin{array}{ccc}
z & \mapsto & \overline{z},\\
w & \mapsto & \overline{w},
\end{array}
$$
which is a lift from $\CC^2\setminus\{x=0,\,u=0\}$ to $M^{\Phi^{\chi}}$ of the following element of $\hbox{Aut}_{CR}(\chi)$:
\begin{equation}
\begin{array}{rrr}
z & \mapsto & z,\\
w & \mapsto & -w.
\end{array}\label{spec1}
\end{equation}
\smallskip\\

$\underline{\hbox{Aut}_{CR}\left(\chi^{(n)}\right)\simeq O_2(\RR)\ltimes\RR^2}:$ This group is generated by the maps
$$
\begin{array}{ccc}
z & \mapsto & \cos\varphi\cdot x+\sin\varphi\cdot u+i\Bigl(\cos(n\varphi)\cdot y+\sin(n\varphi)\cdot v+\beta\Bigr),\\
w & \mapsto & -\sin\varphi\cdot x+\cos\varphi\cdot u+i\Bigl(-\sin(n\varphi)\cdot y+\cos(n\varphi)\cdot v+\gamma\Bigr),
\end{array}
$$
where $\varphi,\beta,\gamma\in\RR$, that form the identity component of $\hbox{Aut}_{CR}\left(\chi^{(n)}\right)$, and map (\ref{spec1}).

\section{The Covers of $\mu_{\alpha}$}\label{mucovers}
\setcounter{equation}{0}

All covers of $\mu_{\alpha}$ can be found by using a map introduced by Rossi in \cite{R}. Let $Q_{+}$ be the variety in $\CC^3$ given by
$$
z_1^2+z_2^2+z_3^2=1.
$$
Consider the map $\Phi^{\mu}: \CC^2\setminus\{0\}\ra Q_{+}$ defined by the formulas
\begin{equation}
\begin{array}{lll}
z_1& = &\displaystyle -i(z^2+w^2)+i\frac{z\overline{w}-w\overline{z}}{|z|^2+|w|^2},\\
\vspace{0cm}&&\\
z_2 & = & \displaystyle z^2-w^2-\frac{z\overline{w}+w\overline{z}}{|z|^2+|w|^2},\\
\vspace{0cm}&&\\
z_3 & = & \displaystyle 2zw+\frac{|z|^2-|w|^2}{|z|^2+|w|^2}.
\end{array}\label{rossioriginal}
\end{equation}
It is straightforward to verify that $\Phi^{\mu}$ is a 2-to-1 covering map onto $Q_{+}\setminus\left(Q_{+}\cap\RR^3\right)$ and that it satisfies
\begin{equation}
\Phi^{\mu}(g(z,w))=\varphi^{\mu}(g)\Phi^{\mu}((z,w)),\label{equivariance}
\end{equation}
for all $g\in SU_2$, $(z,w)\in \CC^2\setminus\{0\}$, where $\varphi^{\mu}$ is the standard 2-to-1 covering homomorphism from $SU_2$ onto $SO_3(\RR)$ defined as follows: for 
$$
g=\left(
\begin{array}{rr}
a & b\\
-\overline{b} & \overline{a}
\end{array}
\right)\in SU_2
$$
(here $|a|^2+|b|^2=1$), set 
\begin{equation}
\varphi^{\mu}(g):=\left(
\begin{array}{rrc}
\hbox{Re}\left(a^2+b^2\right) & \hbox{Im}\left(a^2-b^2\right) & 2\hbox{Im}(ab)\\
-\hbox{Im}\left(a^2+b^2\right) & \hbox{Re}\left(a^2-b^2\right) & 2\hbox{Re}(ab)\\
2\hbox{Im}(a\overline{b}) & -2\hbox{Re}(a\overline{b}) & |a|^2-|b^2|
\end{array}
\right).\label{hommommu}
\end{equation}
In formula (\ref{equivariance}) the actions of $SU_2$ on $\CC^2$ and $SO_3(\RR)$ on $\CC^3$ are standard.

In fact, the covering homomorphism $\varphi^{\mu}$ can be used to give a simple group-theoretic interpretation of Rossi's map $\Phi^{\mu}$.\footnote{This interpretation was suggested to us by one of the referees.}  First of all, we observe that the group $\RR\times SU_2$ acts on $\CC^2\setminus\{0\}$ simply transitively as follows:
$$
(t,g)(z,w):=e^t\cdot g(z,w),
$$
where $t\in\RR$, $g\in SU_2$. On the other hand, the standard action of $SO_3(\RR)$ on $Q_{+}\setminus\left(Q_{+}\cap\RR^3\right)$ can be extended to a simple transitive action of the group $\RR\times SO_3(\RR)$ by diffeomorphisms. Indeed, $\RR\times SO_3(\RR)$ acts simply transitively on
$$
\hat Q_{+}:=\left\{\zeta=(\zeta_1,\zeta_2,\zeta_3)\in\CC^3\setminus\{0\}: \zeta_1^2+\zeta_2^2+\zeta_3^2=0\right\}
$$
as follows:
$$
(t,g)\,\zeta:=e^t\cdot g\,\zeta,
$$
where $t\in\RR$, $g\in SO_3(\RR)$. The manifold $\hat Q_{+}$ is $SO_3(\RR)$-equivariantly diffeomorphic to $Q_{+}\setminus\left(Q_{+}\cap\RR^3\right)$ by means of the map $F_{+}$ given by
$$
\zeta\mapsto \zeta+\xi,
$$
where $\xi\in\RR^3$ is such that $\langle \xi,\xi\rangle_{+}=1$, $\langle \xi,\zeta\rangle_{+}=0$, $\det(\xi,\hbox{Re}\,\zeta,\hbox{Im}\,\zeta)>0$, and $\langle\cdot,\cdot\rangle_{+}$ denotes the standard Hermitian scalar product in $\CC^3$. Using the $SO_3(\RR)$-equivariant diffeomorphism $F_{+}$, we can now push forward the action of $\RR\times SO_3(\RR)$ on $\hat Q_{+}$ to a simple transitive action of $\RR\times SO_3(\RR)$ on $Q_{+}\setminus\left(Q_{+}\cap\RR^3\right)$ by diffeomorphisms.

Thus, as smooth manifolds $\CC^3\setminus\{0\}$ and $Q_{+}\setminus\left(Q_{+}\cap\RR^3\right)$ can be identified with $\RR\times SU_2$ and $\RR\times SO_3(\RR)$, respectively. Then the map $(t,g)\mapsto (t,\varphi^{\mu}(g))$ is precisely Rossi's map $\Phi^{\mu}$ if we choose $(1,0)\in\CC^2\setminus\{0\}$ and $(-i,1,1)\in Q_{+}\setminus\left(Q_{+}\cap\RR^3\right)$ as basepoints. Hence formulas (\ref{rossioriginal}) -- that may look somewhat mysterious at first sight -- are a simple consequence of (\ref{hommommu}).             

Next, we introduce on the domain of $\Phi^{\mu}$ the pull-back complex structure under $\Phi^{\mu}$ and denote the resulting complex manifold by $M^{\Phi^{\mu}}$. This complex structure is invariant under the ordinary action of $SU_2$ on $M^{\Phi^{\mu}}$.
It follows from (\ref{equivariance}) that $\Phi^{\mu}$ maps every $SU_2$-orbit in $M^{\Phi^{\mu}}$ (all such orbits are diffeomorphic to $S^3$) onto an $SO_3(\RR)$-orbit in $Q_{+}\setminus\left(Q_{+}\cap\RR^3\right)$ (we note in passing that $Q_{+}\cap\RR^3$ is also an $SO_3(\RR)$-orbit in $Q_{+}$; it has dimension 2 and does not lie in the range of the map $\Phi^{\mu}$).
Specifically, $\Phi^{\mu}$ maps the $SU_2$-orbit $\left\{(z,w)\in M^{\Phi^{\mu}}: |z|^2+|w|^2=r^2\right\}$, $r>0$, onto the $SO_3(\RR)$-orbit 
$$
\mu_{2r^4+1}^{(2)}:=\left\{(z_1,z_2,z_3)\in\CC^3: |z_1|^2+|z_2|^2+|z_3|^2=2r^4+1 \right\}\cap Q_{+}.
$$
Further, consider a holomorphic map $\Psi^{\mu}:Q_{+}\setminus\left(Q_{+}\cap\RR^3\right)\ra\CC\PP^2\setminus\RR\PP^2$ defined as
$$
(z_1,z_2,z_3)\mapsto (z_1:z_2:z_3).
$$
Clearly, $\Psi^{\mu}$ is a 2-to-1 covering map, and $\Psi^{\mu}\left(\mu_{\alpha}^{(2)}\right)=\mu_{\alpha}$ for every $\alpha>1$. Thus, we have shown that $\tilde\mu_{\alpha}$ coincides with the hypersurface
$$
\mu_{\alpha}^{(4)}:=\left\{(z,w)\in M^{\Phi^{\mu}}: |z|^2+|w|^2=\sqrt{(\alpha-1)/2}\right\},
$$ 
with the $CR$-structure induced from $M^{\Phi^{\mu}}$; the 4-to-1 covering map $\mu_{\alpha}^{(4)}\ra\mu_{\alpha}$ is the composition $\Psi^{\mu}\circ\Phi^{\mu}$.

Next, a straightforward calculation shows that $\Gamma_{\Psi^{\mu}\circ\Phi^{\mu}}$ is a cyclic group of order 4 generated by the map $f^{\mu}$ defined as:
\begin{equation}
\begin{array}{lll}
z & \mapsto & \displaystyle i\frac{z(|z|^2+|w|^2)-\overline{w}}{\sqrt{1+(|z|^2+|w|^2)^2}},\\
\vspace{0cm}&&\\
w & \mapsto & \displaystyle i\frac{w(|z|^2+|w|^2)+\overline{z}}{\sqrt{1+(|z|^2+|w|^2)^2}}.
\end{array}\label{fmu}
\end{equation}
The only non-trivial subgroup of $\Gamma_{\Psi^{\mu}\circ\Phi^{\mu}}$ is then a cyclic subgroup of order 2 generated by $(f^{\mu})^2$. The cover of $\mu_{\alpha}$ corresponding to this subgroup is the hypersurface $\mu_{\alpha}^{(2)}$ with covering map $\mu_{\alpha}^{(2)}\ra\mu_{\alpha}$ coinciding with $\Psi^{\mu}:Q_{+}\setminus\left(Q_{+}\cap\RR^3\right)\ra\CC\PP^2\setminus\RR\PP^2$ and factorization map $\mu_{\alpha}^{(4)}\ra\mu_{\alpha}^{(2)}$ coinciding with $\Phi^{\mu}:M^{\Phi^{\mu}}\ra Q_{+}\setminus\left(Q_{+}\cap\RR^3\right)$.

Thus, every non-trivial (that is, not 1-to-1) cover of $\mu_{\alpha}$ is $CR$-equivalent to either $\mu_{\alpha}^{(4)}$ or $\mu_{\alpha}^{(2)}$. The groups of $CR$-automorphisms of $\mu_{\alpha}^{(4)}$ and $\mu_{\alpha}^{(2)}$ are as follows.
\smallskip\\

$\underline{\hbox{Aut}_{CR}\left(\mu_{\alpha}^{(4)}\right)\simeq SU_2\times_{\hbox{\tiny loc}}\ZZ_4}:$ This group is generated by the maps
$$
\left(
\begin{array}{c}
z\\
w
\end{array}
\right) \mapsto
A
\left(
\begin{array}{c}
z\\
w
\end{array}
\right),
$$
where $A\in SU_2$, that form $\hbox{Aut}_{CR}\left(\mu_{\alpha}^{(4)}\right)^c$, and the map $f^{\mu}$ defined in (\ref{fmu}), which is a lift from $Q_{+}\setminus\RR^3$ to $M^{\Phi^{\mu}}$ of the following element of $\hbox{Aut}_{CR}\left(\mu_{\alpha}^{(2)}\right)^c$: 
\begin{equation}
\begin{array}{rrr}
z_1 & \mapsto & -z_1,\\
z_2 & \mapsto & -z_2,\\
z_3 & \mapsto & -z_3.
\end{array}\label{spec2}
\end{equation}
\smallskip\\

$\underline{\hbox{Aut}_{CR}\left(\mu_{\alpha}^{(2)}\right)\simeq O_3(\RR)}:$  
\begin{equation}
\left(
\begin{array}{c}
z_1\\
z_2\\
z_3
\end{array}
\right)
\mapsto A
\left(
\begin{array}{c}
z_1\\
z_2\\
z_3
\end{array}
\right),\label{autmu2}
\end{equation}
where $A\in O_3(\RR)$.

\section{The Covers of $\nu_{\alpha}$ and $\eta_{\alpha}$}\label{coversnueta}
\setcounter{equation}{0}

In order to find all covers of $\nu_{\alpha}$ and $\eta_{\alpha}$ we introduce an analogue of Rossi's map. Instead of the Hermitian form $|z|^2+|w|^2$ it is associated with the form $|z|^2-|w|^2$. Let $Q_{-}$ be the variety in $\CC^3$ given by
$$
z_1^2+z_2^2-z_3^2=1.
$$
Set $\Omega:=\{(z,w)\in\CC^2:|z|^2-|w|^2\ne 0\}$ and consider the map $\Phi:\Omega\ra Q_{-}$ defined by the formulas
\begin{equation}
\begin{array}{lll}
z_1 & = &\displaystyle -i(z^2+w^2)-i\frac{z\overline{w}+w\overline{z}}{|z|^2-|w|^2},\\
\vspace{0cm}&&\\
z_2 & = & \displaystyle z^2-w^2+\frac{z\overline{w}-w\overline{z}}{|z|^2-|w|^2},\\
\vspace{0cm}&&\\
z_3 & = &\displaystyle -2izw-i\frac{|z|^2+|w|^2}{|z|^2-|w|^2}.
\end{array}\label{analoguerossi}
\end{equation}
It is straightforward to verify that the range of $\Phi$ is $Q_{-}\setminus\left(Q_{-}\cap{\cal W}\right)$, where
$$
\begin{array}{lll}
{\cal W}:&=&i\RR^3\cup \RR^3\cup\\
&&\displaystyle\left\{(z_1,z_2,z_3)\in\CC^3\setminus\RR^3:|iz_1+z_2|=|iz_3-1|,\,|iz_1-z_2|=|iz_3+1|\right\},
\end{array}
$$
and that the restrictions of $\Phi$ to the domains 
$$
\begin{array}{lll}
\Omega^{>} & := & \{(z,w)\in\CC^2:|z|^2-|w|^2> 0\},\\
\Omega^{<} & := & \{(z,w)\in\CC^2:|z|^2-|w|^2< 0\}
\end{array}
$$
are 2-to-1 covering maps onto $\Phi(\Omega^{>})$ and $\Phi(\Omega^{<})$, respectively (note that $\Phi(\Omega^{<})$ is obtained from $\Phi(\Omega^{>})$ by applying the transformation $z_1\mapsto -z_1$, $z_2\mapsto z_2$, $z_3\mapsto -z_3$).

The map $\Phi$ satisfies
\begin{equation}
\Phi(g(z,w))=\varphi(g)\Phi((z,w)),\label{equivariance1}
\end{equation}
for all $g\in SU_{1,1}$, $(z,w)\in\Omega$, where $\varphi$ is the standard 2-to-1 covering homomorphism from $SU_{1,1}$ onto $SO_{2,1}(\RR)^c$, defined as follows: for 
$$
g=\left(
\begin{array}{rr}
a & b\\
\overline{b} & \overline{a}
\end{array}
\right)\in SU_{1,1}
$$
(here $|a|^2-|b|^2=1$), set 
\begin{equation}
\varphi(g):=\left(
\begin{array}{rrc}
\hbox{Re}\left(a^2+b^2\right) & \hbox{Im}\left(a^2-b^2\right) & 2\hbox{Re}(ab)\\
-\hbox{Im}\left(a^2+b^2\right) & \hbox{Re}\left(a^2-b^2\right) & -2\hbox{Im}(ab)\\
2\hbox{Re}(a\overline{b}) & 2\hbox{Im}(a\overline{b}) & |a|^2+|b^2|
\end{array}
\right).\label{homomphi}
\end{equation}
The actions of $SU_{1,1}$ on $\CC^2$ and  $SO_{2,1}(\RR)^c$ on $\CC^3$ in formula (\ref{equivariance1}) are standard.

Hence $\Phi$ maps every $SU_{1,1}$-orbit in $\Omega$ onto an $SO_{2,1}(\RR)^c$-orbit in $Q_{-}\setminus\left(Q_{-}\cap{\cal W}\right)$. Note that $SO_{2,1}(\RR)^c$ has exactly four orbits in $Q_{-}$ that do not lie in the range of the map $\Phi$:
$$
\begin{array}{lll}
O_1 &:=& i\RR^3\cap Q_{-}=\left\{(z_1,z_2,z_3)\in\CC^3: |z_1|^2+|z_2|^2-|z_3|^2=-1 \right\}\cap Q_{-},\\
O_2 &:=& \RR^3\cap Q_{-},\\
O_3 &:=& \Bigl\{(z_1,z_2,z_3)\in\CC^3\setminus\RR^3:|iz_1+z_2|=|iz_3-1|,\\
&&|iz_1-z_2|=|iz_3+1|,\,\hbox{Im}\,z_3<0\Bigr\}\cap Q_{-},\\
O_4 &:=& \Bigl\{(z_1,z_2,z_3)\in\CC^3\setminus\RR^3:|iz_1+z_2|=|iz_3-1|,\\
&&|iz_1-z_2|=|iz_3+1|,\,\hbox{Im}\,z_3>0\Bigr\}\cap Q_{-}.
\end{array}
$$
The orbits $O_1$, $O_2$ are 2-dimensional, the orbits $O_3$, $O_4$ are 3-dimensional, and $O_2$, $O_3$, $O_4$ lie in the set
$$
{\frak S}:=\left\{(z_1,z_2,z_3)\in\CC^3: |z_1|^2+|z_2|^2-|z_3|^2=1 \right\}\cap Q_{-}.
$$
In fact, we have
$$
{\frak S}=O_2\cup O_3\cup O_4\cup O_5\cup O_6,
$$
where $O_5$, $O_6$ are the following 3-dimensional $SO_{2,1}(\RR)^c$-orbits in $Q_{-}$:
$$
\begin{array}{lll}
O_5&:=&\Bigl\{(z_1,z_2,z_3)\in\CC^3\setminus\RR^3:|iz_1+z_2|=|iz_3+1|,\\
&&|iz_1-z_2|=|iz_3-1|,\,\hbox{Im}\,z_3<0\Bigr\}\cap Q_{-},\\
O_6&:=&\Bigl\{(z_1,z_2,z_3)\in\CC^3\setminus\RR^3:|iz_1+z_2|=|iz_3+1|,\\
&&|iz_1-z_2|=|iz_3-1|,\,\hbox{Im}\,z_3>0\Bigr\}\cap Q_{-}.
\end{array}
$$
In contrast with $O_2$, $O_3$, $O_4$, however, the orbits $O_5$, $O_6$ lie in the range of $\Phi$ and are the images under $\Phi$ of the sets $\{(z,w)\in\CC^2: |z|^2-|w|^2=1\}$, $\{(z,w)\in\CC^2: |z|^2-|w|^2=-1\}$, respectively.

From now on we will only consider the restriction $\Phi^{>}$ of $\Phi$ to $\Omega^{>}$. The range of $\Phi^{>}$ is $D^{>}:=\Sigma^{\nu}\cup\Sigma^{\eta}\cup O_5$, where
$$
\begin{array}{lll}
\Sigma^{\nu} & := & \Bigl\{(z_1,z_2,z_3)\in\CC^3: -1<|z_1|^2+|z_2|^2-|z_3|^2<1,\\
&&\hbox{Im}\,z_3<0\Bigr\}\cap Q_{-},
\end{array}
$$
and
$$
\begin{array}{lll}
\Sigma^{\eta}& := & \Bigl\{(z_1,z_2,z_3)\in\CC^3:|z_1|^2+|z_2|^2-|z_3|^2>1,\\
&&\hbox{Im}(z_2(\overline{z_1}+\overline{z_3}))>0\Bigr\}\cap Q_{-}.
\end{array}
$$

The covering homomorphism $\varphi$ defined in (\ref{homomphi}) can be used to give a group-theoretic interpretation of the map $\Phi^{>}$ analogous to that of Rossi's map $\Phi^{\mu}$ from the previous section, where the homomorphism $\varphi^{\mu}$ defined in (\ref{hommommu}) was utilized. The group $\RR\times SU_{1,1}$ acts on $\Omega^{>}$ simply transitively as follows:
$$
(t,g)(z,w):=e^t\cdot g(z,w),
$$
where $t\in\RR$, $g\in SU_{1,1}$. On the other hand, the standard action of $SO_{2,1}(\RR)^c$ on $D^{>}$ can be extended to a simple transitive action of the group $\RR\times SO_{2,1}(\RR)^c$ by diffeomorphisms. Indeed, $\RR\times SO_{2,1}(\RR)^c$ acts simply transitively on each of the two connected components of the set
$$
\hat Q_{-}:=\left\{\zeta=(\zeta_1,\zeta_2,\zeta_3)\in\CC^3: \zeta_1^2+\zeta_2^2+\zeta_3^2=0,\,\langle\zeta,\zeta\rangle_{-}>0\right\},
$$
where $\langle\zeta,\zeta'\rangle_{-}:=\zeta_1\overline{\zeta_1'}+\zeta_2\overline{\zeta_2'}-\zeta_3\overline{\zeta_3'}$. 
The action is given as follows:
$$
(t,g)\,\zeta:=e^t\cdot g\,\zeta,
$$
where $t\in\RR$, $g\in SO_{2,1}(\RR)^c$. Let $\hat Q_{-}^0$ be the connected component of $\hat Q_{-}$ that contains the point $(-i,1,0)$.   
The manifold $\hat Q_{-}^0$ is $SO_{2,1}(\RR)^c$-equivariantly diffeomorphic to $D^{>}$ by means of the map $F_{-}$ defined as
$$
\zeta\mapsto \zeta+i\xi,
$$
where $\xi\in\RR^3$ is such that $\langle \xi,\xi\rangle_{-}=-1$, $\langle \xi,\zeta\rangle_{-}=0$, $\det(\xi,\hbox{Re}\,\zeta,\hbox{Im}\,\zeta)<0$. Using the $SO_{2,1}(\RR)^c$-equivariant diffeomorphism $F_{-}$, we can now push forward the action of $\RR\times SO_{2,1}(\RR)^c$ on $\hat Q_{-}^0$ to a simple transitive action of $\RR\times SO_{2,1}(\RR)^c$ on $D^{>}$ by diffeomorphisms.

Thus, as smooth manifolds $\Omega^{>}$ and $D^{>}$ can be identified with $\RR\times SU_{1,1}$ and $\RR\times SO_{2,1}(\RR)^c$, respectively. Then the map $(t,g)\mapsto (t,\varphi(g))$ is exactly the map $\Phi^{>}$ if we choose $(1,0)\in\Omega^{>}$ and $(-i,1,-i)\in D^{>}$ as basepoints. 

We will now concentrate on the domains in $\Omega^{>}$ lying above $Q_{-}\setminus(O_1\cup{\frak S})$. Let $\Phi^{\nu}$, $\Phi^{\eta}$ denote the restrictions of $\Phi^{>}$ to the domains
$$
\Omega^{\nu}:=\{(z,w)\in\CC^2:0<|z|^2-|w|^2<1\}
$$
and
$$
\Omega^{\eta}:=\{(z,w)\in\CC^2:|z|^2-|w|^2>1\},
$$
respectively. The maps $\Phi^{\nu}$ and $\Phi^{\eta}$ are 2-to-1 covering maps onto $\Sigma^{\nu}$ and $\Sigma^{\eta}$. Introduce on $\Omega^{\nu}$, $\Omega^{\eta}$ the pull-back complex structures under the maps $\Phi^{\nu}$, $\Phi^{\eta}$, respectively, and denote the resulting complex manifolds by $M^{\Phi^{\nu}}$, $M^{\Phi^{\eta}}$. These complex structures are invariant under the ordinary action of $SU_{1,1}$. The map $\Phi^{\nu}$ takes the $SU_{1,1}$-orbit
\begin{equation}
\nu_{2r^4-1}^{(2)}:=\left\{(z,w)\in M^{\Phi^{\nu}}: |z|^2-|w|^2=r^2\right\}\label{nu2}
\end{equation}
onto the $SO_{2,1}(\RR)^c$-orbit 
$$
\nu'_{2r^4-1}:=\left\{(z_1,z_2,z_3)\in\Sigma^{\nu}: |z_1|^2+|z_2|^2-|z_3|^2=2r^4-1\right\},
$$
where $0<r<1$. Similarly, the map $\Phi^{\eta}$ takes the $SU_{1,1}$-orbit
\begin{equation}
\eta_{2r^4-1}^{(4)}:=\left\{(z,w)\in M^{\Phi^{\eta}}: |z|^2-|w|^2=r^2\right\}\label{eta4}
\end{equation}
onto the $SO_{2,1}(\RR)^c$-orbit 
\begin{equation}
\eta_{2r^4-1}^{(2)}:=\left\{(z_1,z_2,z_3)\in\Sigma^{\eta}: |z_1|^2+|z_2|^2-|z_3|^2=2r^4-1\right\},\label{eta2}
\end{equation}
where $r>1$. 

Observe now that $z_3\ne 0$ on $\Sigma^{\nu}$ and consider the holomorphic map $\Psi^{\nu}:\Sigma^{\nu}\ra\CC^2$ defined as follows:
$$
\begin{array}{lll}
z & = & \displaystyle z_1/z_3,\\
w & = & \displaystyle z_2/z_3.
\end{array}
$$
This map is 1-to-1, it takes $\Sigma^{\nu}$ onto
$$
D^{\nu}:=\Bigl\{(z,w)\in\CC^2: -|z^2+w^2-1|<|z|^2+|w|^2-1<|z^2+w^2-1|\Bigr\},
$$
and establishes equivalence between $\nu'_{\alpha}$ and $\nu_{\alpha}$ for $-1<\alpha<1$. Next, we note that $z_1\ne 0$ on $\Sigma^{\eta}$ and consider the holomorphic map $\Psi^{\eta}:\Sigma^{\eta}\ra\CC^2$ defined by
$$
\begin{array}{lll}
z & = & \displaystyle z_2/z_1,\\
w & = & \displaystyle z_3/z_1.
\end{array}
$$
It is easy to see that $\Psi^{\eta}$ is a 2-to-1 covering map onto 
$$
D^{\eta}:=\Bigl\{(z,w)\in\CC^2: 1+|z|^2-|w|^2>|1+z^2-w^2|,\,\hbox{Im}(z(1+\overline{w}))>0\Bigr\},
$$
and realizes $\eta_{\alpha}^{(2)}$ as a 2-sheeted cover of $\eta_{\alpha}$ for $\alpha>1$.

Let $\Lambda:\CC\times\Delta\to\Omega^{>}$ be the following covering map:
$$
\begin{array}{lll}
z & := & e^s,\\
w & := & e^st,
\end{array}
$$
where $s\in\CC$, $t\in\Delta$ and $\Delta$ is the unit disk. Further, define
$$
\begin{array}{lll}
U^{\nu} & := & \left\{(s,t)\in\CC^2:|t|<1,\,\exp(2\hbox{Re}\,s)(1-|t|^2)<1\right\},\\
U^{\eta} & := & \left\{(s,t)\in\CC^2:|t|<1,\,\exp(2\hbox{Re}\,s)(1-|t|^2)>1\right\}.
\end{array}
$$
Denote by $\Lambda^{\nu}$, $\Lambda^{\eta}$ the restrictions of $\Lambda$ to $U^{\nu}$, $U^{\eta}$, respectively. Clearly, $U^{\nu}$ covers $M^{\Phi^{\nu}}$ by means of $\Lambda^{\nu}$, and $U^{\eta}$ covers $M^{\Phi^{\eta}}$ by means of $\Lambda^{\eta}$. Introduce now on $U^{\nu}$, $U^{\eta}$ the pull-back complex structures under the maps $\Lambda^{\nu}$, $\Lambda^{\eta}$, respectively, and denote the resulting complex manifolds by $M^{\Lambda^{\nu}}$, $M^{\Lambda^{\eta}}$. Then the simply-connected hypersurface
$$
\left\{(s,t)\in M^{\Lambda^{\nu}}:r^2\exp\left(-2\hbox{Re}\,s\right)+|t|^2=1\right\}
$$
covers $\nu_{2r^4-1}^{(2)}$ by means of the map $\Lambda^{\nu}$ for $0<r<1$, and the simply-connected hypersurface
$$
\left\{(s,t)\in M^{\Lambda^{\eta}}:r^2\exp\left(-2\hbox{Re}\,s\right)+|t|^2=1\right\}
$$
covers $\eta_{2r^4-1}^{(4)}$ by means of the map $\Lambda^{\eta}$ for $r>1$.

Thus, $\tilde\nu_{\alpha}$ for $-1<\alpha<1$ coincides with
\begin{equation}
\nu_{\alpha}^{(\infty)}:=\left\{(s,t)\in M^{\Lambda^{\nu}}:\sqrt{(\alpha+1)/2}\exp\left(-2\hbox{Re}\,s\right)+|t|^2=1\right\},\label{universalnu}
\end{equation}
with the $CR$-structure induced from the complex structure of $M^{\Lambda^{\nu}}$, and $\tilde\eta_{\alpha}$ for $\alpha>1$ coincides with
\begin{equation}
\eta_{\alpha}^{(\infty)}:=\left\{(s,t)\in M^{\Lambda^{\eta}}:\sqrt{(\alpha+1)/2}\exp\left(-2\hbox{Re}\,s\right)+|t|^2=1\right\},\label{universaleta}
\end{equation}
with the $CR$-structure induced from the complex structure of $M^{\Lambda^{\eta}}$. The covering maps $\nu_{\alpha}^{(\infty)}\ra\nu_{\alpha}$ and $\eta_{\alpha}^{(\infty)}\ra\eta_{\alpha}$ are respectively $\Psi^{\nu}\circ\Phi^{\nu}\circ\Lambda^{\nu}:U^{\nu}\ra D^{\nu}$ and $\Psi^{\eta}\circ\Phi^{\eta}\circ\Lambda^{\eta}:U^{\eta}\ra D^{\eta}$.   

Next, the group $\Gamma_{\Psi^{\nu}\circ\Phi^{\nu}\circ\Lambda^{\nu}}=\Gamma_{\Phi^{\nu}\circ\Lambda^{\nu}}$ consists of the maps
\begin{equation}
\begin{array}{lll}
s & \mapsto & s+\pi i k,\quad k\in\ZZ,\\
t & \mapsto & t.
\end{array}\label{covmaps}
\end{equation}
Let $\Gamma\subset\Gamma_{\Phi^{\nu}\circ\Lambda^{\nu}}$ be a subgroup. Then there exists an integer $n \ge 0$ such that every element of $\Gamma$ has the form
\begin{equation}
\begin{array}{lll}
s & \mapsto & s+\pi i n k,\quad k\in\ZZ,\\
t & \mapsto & t.
\end{array}\label{pin}
\end{equation}
Suppose that $n\ge 2$ and set
$$
\Omega^{\nu\,(n)}:=\left\{(z,w)\in\CC^2: 0<|z|^n-|z|^{n-2}|w|^2<1\right\}
$$
(note that $\Omega^{\nu\,(2)}=\Omega^{\nu}$). Consider the map $\Phi^{\nu}_n$ from $\Omega^{\nu\,(n)}$ to $\Sigma^{\nu}$ defined as follows:
\begin{equation}
\begin{array}{lll}
z_1 & = &\displaystyle -i(z^n+z^{n-2}w^2)-i\frac{z\overline{w}+w\overline{z}}{|z|^2-|w|^2},\\
\vspace{0cm}&&\\
z_2 & = & \displaystyle z^n-z^{n-2}w^2+\frac{z\overline{w}-w\overline{z}}{|z|^2-|w|^2},\\
\vspace{0cm}&&\\
z_3 & = &\displaystyle -2iz^{n-1}w-i\frac{|z|^2+|w|^2}{|z|^2-|w|^2}.
\end{array}\label{rossin}
\end{equation}
This map is a generalization of map (\ref{analoguerossi}) introduced at the beginning of the section and also can be viewed as an analogue of Rossi's map (\ref{rossioriginal}). The extension of this map by the same formula to all of $\Omega^{>}$ admits a group-theoretic interpretation analogous to those given above for the maps $\Phi^{\mu}$ and $\Phi^{>}$. It uses an $n$-to-1 covering homomorphism $SO_{2,1}(\RR)^{c\,(n)}\ra SO_{2,1}(\RR)^c$, where the $n$-sheeted cover $SO_{2,1}(\RR)^{c\,(n)}$ of $SO_{2,1}(\RR)^c$ is realized as the group of maps of the form (\ref{group2}) that will appear below, acting on $\Omega^{>}$ (note that this group reduces to the group $SU_{1,1}$ for $n=2$). We do not provide a detailed construction here since it is very similar to that for the map $\Phi^{>}$.

Denote by $M^{\Phi^{\nu}_n}$ the domain $\Omega^{\nu\,(n)}$ with the pull-back complex structure under the map $\Phi^{\nu}_n$ (note that $M^{\Phi^{\nu}_2}=M^{\Phi^{\nu}}$). Then the hypersurface 
\begin{equation}
\nu_{\alpha}^{(n)}:=\left\{(z,w)\in M^{\Phi^{\nu}_n}:|z|^n-|z|^{n-2}|w|^2=\sqrt{(\alpha+1)/2}\right\},\label{finitenu}
\end{equation}
equipped with the $CR$-structure induced by the complex structure of $M^{\Phi^{\nu}_n}$ is an $n$-sheeted cover of $\nu_{\alpha}$ corresponding to $\Gamma$ with covering map $\nu_{\alpha}^{(n)}\ra\nu_{\alpha}$ coinciding with $\Psi^{\nu}\circ\Phi^{\nu}_n:M^{\Phi^{\nu}_n}\ra D^{\nu}$ and factorization map $\nu_{\alpha}^{(\infty)}\ra\nu_{\alpha}^{(n)}$ given by
\begin{equation}
\begin{array}{lll}
z & \mapsto & e^{2s/n},\\
w & \mapsto & e^{2s/n}t,
\end{array}\label{expcover}
\end{equation}
(observe that for $n=2$ formula (\ref{finitenu}) coincides with (\ref{nu2})). Thus, every non-trivial cover of $\nu_{\alpha}$ is $CR$-equivalent to either $\nu_{\alpha}^{(\infty)}$ or $\nu_{\alpha}^{(n)}$ for some $n\in\NN$, $n\ge 2$. The groups of $CR$-automorphisms of $\nu_{\alpha}^{(\infty)}$ and $\nu_{\alpha}^{(n)}$ are given below.
\smallskip\\

\noindent {\bf (A)} $\underline{\hbox{Aut}_{CR}\left(\nu_{\alpha}^{(\infty)}\right)\simeq\widetilde{SO}_{2,1}(\RR)^c\rtimes_{\hbox{\tiny loc}}\ZZ}:$ This group is generated by the following maps (they form the subgroup $\hbox{Aut}_{CR}\left(\nu_{\alpha}^{(\infty)}\right)^c$):
\begin{equation}
\begin{array}{lll}
s & \mapsto & \displaystyle s+\ln(a+bt),\\
\vspace{0cm}&&\\
t & \mapsto & \displaystyle\frac{\overline{b}+\overline{a}t}{a+bt},
\end{array}\label{group1}
\end{equation}
where $\ln$ is any branch of the logarithm and $|a|^2-|b|^2=1$, and the map
\begin{equation}
\begin{array}{lll}
s & \mapsto & \displaystyle \overline{s} +\ln'\left(-\frac{1+e^{2s}t(1-|t|^2)}{\sqrt{1-\exp\left(4\hbox{Re}\,s\right)(1-|t|^2)^2}}\right),\\
\vspace{0cm}&&\\
t & \mapsto & \displaystyle -\frac{\overline{t}+e^{2s}(1-|t|^2)}{1+e^{2s}t(1-|t|^2)},
\end{array}\label{spec4}
\end{equation}
for some branch $\ln'$ of the logarithm. Map (\ref{spec4}) is a lift from $\Sigma^{\nu}$ to $M^{\Lambda^{\nu}}$ of the following element of $\hbox{Aut}_{CR}(\nu'_{\alpha})$:
\begin{equation}
\begin{array}{rrr}
z_1 & \mapsto & -z_1,\\
z_2 & \mapsto & z_2,\\
z_3 & \mapsto & z_3.
\end{array}\label{spec3}
\end{equation}
Here $\rtimes_{\hbox{\tiny loc}}$ denotes local semidirect product.
\smallskip\\

\noindent {\bf (B)} $\underline{\hbox{Aut}_{CR}\left(\nu_{\alpha}^{(n)}\right)\simeq SO_{2,1}(\RR)^{c\,(n)}\rtimes_{\hbox{\tiny loc}}\ZZ_{2n}}:$ This group is generated by the following maps (they form the subgroup $\hbox{Aut}_{CR}\left(\nu_{\alpha}^{(n)}\right)^c$):  
\begin{equation}
\begin{array}{lll}
z & \mapsto & \displaystyle z\,\sqrt[n]{\left(a+b\, w/z\right)^2}, \\
\vspace{0cm}&&\\ 
w & \mapsto & \displaystyle z\,\frac{\overline{b}+\overline{a}\,w/z}{a+b\,w/z}\sqrt[n]{\left(a+b\,w/z\right)^2},
\end{array}\label{group2}
\end{equation}
where $\sqrt[n]{}$ is any branch of the $n$th root and $|a|^2-|b|^2=1$, and the map
\begin{equation}
\begin{array}{lll}
z & \mapsto & \displaystyle\overline{z}\left(\sqrt[n]{\frac{\Bigl(1+z^{n-1}w(1-|w|^2/|z|^2)\Bigr)^2}{1-|z|^{2n}(1-|w|^2/|z|^2)^2}}\right)',\\
\vspace{0cm}&&\\
w & \mapsto & -\displaystyle\frac{\overline{w}/\overline{z}+z^n(1-|w|^2/|z|^2)}{1+z^{n-1}w(1-|w|^2/|z|^2)}\times\\
\vspace{0cm}&&\\
&&\displaystyle\hspace{1.2cm}\overline{z}\left(\sqrt[n]{\frac{\Bigl(1+z^{n-1}w(1-|w|^2/|z|^2)\Bigr)^2}{1-|z|^{2n}(1-|w|^2/|z|^2)^2}}\right)',
\end{array}\label{spec5}
\end{equation}
for some branch $\left(\sqrt[n]{}\right)'$ of $\sqrt[n]{}$\,. Map (\ref{spec5}) is a lift from $\Sigma^{\nu}$ to $M^{\Phi^{\nu}_n}$ of map (\ref{spec3}) (recall that $SO_{2,1}^{c\,(n)}(\RR)$ is the $n$-sheeted cover of $SO_{2,1}(\RR)^c$).
\smallskip\\ 

Further, the group $\Gamma_{\Psi^{\eta}\circ\Phi^{\eta}\circ\Lambda^{\eta}}$ is generated by all maps of the form (\ref{covmaps}) and the map $f^{\eta}$ defined as follows: 
\begin{equation}
\begin{array}{lll}
s & \mapsto & \displaystyle 2s+\overline{s} +\ln'\left(i\frac{1-|t|^2+e^{-2s}\overline{t}}{\sqrt{\exp\left(4\hbox{Re}\,s\right)(1-|t|^2)^2-1}}\right),\\
\vspace{0cm}&&\\
t & \mapsto & \displaystyle \frac{1+e^{2s}t(1-|t|^2)}{\overline{t}+e^{2s}(1-|t|^2)},
\end{array}\label{spec6}
\end{equation}
for some fixed branch $\ln'$ of the logarithm. The map $f^{\eta}$ is a lift from $\Sigma^{\eta}$ to $M^{\Lambda^{\eta}}$ of the element of $\hbox{Aut}_{CR}\left(\eta_{\alpha}^{(2)}\right)$ given by formula (\ref{spec2}). At the same time, $f^{\eta}$ is a lift from $M^{\Phi^{\eta}}$ to $M^{\Lambda^{\eta}}$ of the map
\begin{equation}
\begin{array}{lll}
z & \mapsto & \displaystyle i\frac{z(|z|^2-|w|^2)+\overline{w}}{\sqrt{(|z|^2-|w|^2)^2-1}},\\
\vspace{0cm}&&\\
w & \mapsto & \displaystyle i\frac{w(|z|^2-|w|^2)+\overline{z}}{\sqrt{(|z|^2+|w|^2)^2-1}}.
\end{array}\label{specc}
\end{equation}
Since the square of map (\ref{specc}) is
\begin{equation}
\begin{array}{lll}
z & \mapsto & -z,\\
w & \mapsto & -w,
\end{array}\label{special}
\end{equation}
it follows that $(f^{\eta})^2$ is a lift from $M^{\Phi^{\eta}}$ to $M^{\Lambda^{\eta}}$ of map (\ref{special}) and thus has the form   
(\ref{covmaps}) with $k$ odd. Since $f^{\eta}$ clearly commutes with all maps (\ref{covmaps}), the group $\Gamma_{\Psi^{\eta}\circ\Phi^{\eta}\circ\Lambda^{\eta}}$ is isomorphic to $\ZZ$ and is generated by $f^{\eta}\circ g$, where $g$ has the form (\ref{covmaps}).

Let $\Gamma\subset\Gamma_{\Psi^{\eta}\circ\Phi^{\eta}\circ\Lambda^{\eta}}$ be a subgroup. It then follows that $\Gamma$ is generated by either a map of the form (\ref{covmaps}) or by $f^{\eta}\circ h$, where $h$ has the form (\ref{covmaps}). In the first case there exists an integer $n\ge 0$ such that every element of $\Gamma$ has the form (\ref{pin}). Suppose that $n\ge 2$ and set
$$
\Omega^{\eta\,(n)}:=\left\{(z,w)\in\CC^2: |z|^n-|z|^{n-2}|w|^2>1\right\}
$$
(note that $\Omega^{\eta\,(2)}=\Omega^{\eta}$). Consider the map $\Phi^{\eta}_n$ from $\Omega^{\eta\,(n)}$ to $\Sigma^{\eta}$ defined by formula (\ref{rossin}). Denote by $M^{\Phi^{\eta}_n}$ the domain $\Omega^{\eta\,(n)}$ with the pull-back complex structure under the map $\Phi^{\eta}_n$ (note that $M^{\Phi^{\eta}_2}=M^{\Phi^{\eta}}$). Then the hypersurface 
\begin{equation}
\eta_{\alpha}^{(2n)}:=\left\{(z,w)\in M^{\Phi^{\eta}_n}:|z|^n-|z|^{n-2}|w|^2=\sqrt{(\alpha+1)/2}\right\},\label{finiteeta1}
\end{equation}
equipped with the $CR$-structure induced by the complex structure of $M^{\Phi^{\eta}_n}$ is a $2n$-sheeted cover of $\eta_{\alpha}$ corresponding to $\Gamma$ with covering map $\eta_{\alpha}^{(2n)}\ra\eta_{\alpha}$ coinciding with $\Psi^{\eta}\circ\Phi^{\eta}_n:M^{\Phi^{\eta}_n}\ra D^{\eta}$ and factorization map $\eta_{\alpha}^{(\infty)}\ra\eta_{\alpha}^{(2n)}$ given by formula (\ref{expcover}); note that for $n=2$ formula (\ref{finiteeta1}) coincides with (\ref{eta4}). For $n=1$ we obtain the hypersurface $\eta_{\alpha}^{(2)}$ defined in (\ref{eta2}) that covers $\eta_{\alpha}$ by means of the 2-to-1 map $\Psi^{\eta}:\Sigma^{\eta}\ra D^{\eta}$; the factorization map $\eta_{\alpha}^{(\infty)}\ra\eta_{\alpha}^{(2)}$ is $\Phi^{\eta}\circ\Lambda^{\eta}:M^{\Lambda^{\eta}}\ra\Sigma^{\eta}$. 

We will now describe finitely-sheeted covers of $\eta_{\alpha}$ of odd orders. They arise in the case when $\Gamma$ is generated by $f^{\eta}\circ h$, where $h$ has the form (\ref{covmaps}). In this case the group $\Gamma$ can be represented as $\Gamma=\Gamma_0\cup\Bigl((f^{\eta}\circ h)\circ\Gamma_0\Bigr)$, where $\Gamma_0$ is a subgroup that consists of maps of the form (\ref{pin}) for an odd positive integer $n$. We will first factor $M^{\Phi^{\eta}}$ with respect to the action of the subgroup of $\Gamma_0$ corresponding to even $k$. Namely,  $M^{\Phi^{\eta}}$ by means of the map
\begin{equation}
\begin{array}{lll}
z & \mapsto & e^{s/n},\\
w & \mapsto & e^{s/n}t,
\end{array}\label{expmm}
\end{equation}
covers the manifold $M^{\Phi^{\eta}_{2n}}$ and, accordingly, $\eta_{\alpha}^{(\infty)}$ covers $\eta_{\alpha}^{(4n)}$. In order to obtain the cover of $\eta_{\alpha}$ corresponding to the group $\Gamma$, the hypersurface $\eta_{\alpha}^{(4n)}$ must be further factored by the action of the cyclic group of four elements generated by the following automorphism $f^{\eta}_n$ of $M^{\Phi^{\eta}_{2n}}$:
$$
\begin{array}{lll}
z & \mapsto & \displaystyle iz^2\overline{z}\left(\sqrt[n]{\frac{1-|w|^2/|z|^2+z^{-2n}\overline{w}/\overline{z}}{\sqrt{|z|^{4n}(1-|w|^2/|z|^2)^2-1}}}\right)',\\
\vspace{0cm}&&\\
w & \mapsto & i\displaystyle\frac{1+z^{2n-1}w(1-|w|^2/|z|^2)}{\overline{w}/\overline{z}+z^{2n}(1-|w|^2/|z|^2)}\times\\
\vspace{0cm}&&\\
&&\displaystyle z^2\overline{z}\left(\sqrt[n]{\frac{1-|w|^2/|z|^2+z^{-2n}\overline{w}/\overline{z}}{\sqrt{|z|^{4n}(1-|w|^2/|z|^2)^2-1}}}\right)',
\end{array}
$$
for some branch $\left(\sqrt[n]{}\right)'$ of $\sqrt[n]{}$\,. Let $\hat M^{\Phi^{\eta}_n}$  denote the manifold arising from $M^{\Phi^{\eta}_{2n}}$ by means of this factorization and $\Pi^{\eta}_n: M^{\Phi^{\eta}_{2n}}\ra\hat M^{\Phi^{\eta}_n}$ denote the corresponding 4-to-1 factorization map. Next, define
\begin{equation}
\eta_{\alpha}^{(n)}:=\Pi^{\eta}_n\left(\eta_{\alpha}^{(4n)}\right).\label{etanodd}
\end{equation}
The hypersurface $\eta_{\alpha}^{(n)}$ with the $CR$-structure induced from the complex structure of $\hat M^{\Phi^{\eta}_n}$ is an $n$-sheeted cover of $\eta_{\alpha}$ corresponding to $\Gamma$ with the factorization map $\eta_{\alpha}^{(\infty)}\ra\eta_{\alpha}^{(n)}$ coinciding with the composition of map (\ref{expmm}) and $\Pi^{\eta}_n$. Note that for $n=1$ the map $f^{\eta}_n=f^{\eta}_1$ coincides with the map defined in (\ref{specc}), and $\Pi^{\eta}_n=\Pi^{\eta}_1$ coincides with $\Psi^{\eta}\circ\Phi^{\eta}$. Both $\Pi^{\eta}_n$ and the covering map $\eta_{\alpha}^{(n)}\ra\eta_{\alpha}$ (which extends to a covering map $\hat M^{\Phi^{\eta}_n}\ra D^{\eta}$) can be computed explicitly for any odd $n\in\NN$, but, since the resulting formulas are quite lengthly and not very instructive, we omit them.

We will now write down the groups of $CR$-automorphisms of $\eta_{\alpha}^{(\infty)}$ and $\eta_{\alpha}^{(n)}$.
\smallskip\\

\noindent {\bf (C)} $\underline{\hbox{Aut}_{CR}\left(\eta_{\alpha}^{(\infty)}\right)\simeq \widetilde{SO}_{2,1}(\RR)^c\times_{\hbox{\tiny loc}}\ZZ}:$ This group is generated by its connected identity component that consists of all maps of the form (\ref{group1}), and the map $f^{\eta}$ defined in (\ref{spec6}).
\smallskip\\

\noindent {\bf (D)} $\underline{\hbox{Aut}_{CR}\left(\eta_{\alpha}^{(2)}\right)\simeq SO_{2,1}(\RR)^c\times\ZZ_2}:$ This group is generated by its connected identity component that consists of all maps of the form (\ref{autmu2}), where $A\in SO_{2,1}(\RR)^c$ and the map given by formula (\ref{spec2}).
\smallskip\\

\noindent {\bf (D')} $\underline{\hbox{Aut}_{CR}\left(\eta_{\alpha}^{(2n)}\right)\simeq SO_{2,1}(\RR)^{c\,(n)}\times_{\hbox{\tiny loc}}\ZZ_{4n},\, n\ge 2}:$ This group is generated by its connected identity component that consists of all maps of the form (\ref{group2}), and the map
$$
\begin{array}{lll}
z & \mapsto & \displaystyle z^2\overline{z}\left(\sqrt[n]{\frac{\Bigl(1-|w|^2/|z|^2+z^{-n}\overline{w}/\overline{z}\Bigr)^2}{|z|^{2n}(1-|w|^2/|z|^2)^2-1}}\right)',\\
\vspace{0cm}&&\\
w & \mapsto & \displaystyle\frac{1+z^{n-1}w(1-|w|^2/|z|^2)}{\overline{w}/\overline{z}+z^{n}(1-|w|^2/|z|^2)}\times\\
\vspace{0cm}&&\\
&&\displaystyle z^2\overline{z}\left(\sqrt[n]{\frac{\Bigl(1-|w|^2/|z|^2+z^{-n}\overline{w}/\overline{z}\Bigr)^2}{|z|^{2n}(1-|w|^2/|z|^2)^2-1}}\right)',
\end{array}
$$
for some branch $\left(\sqrt[n]{}\right)'$ of $\sqrt[n]{}$.
\smallskip\\

\noindent {\bf (E)} $\underline{\hbox{Aut}_{CR}\left(\eta_{\alpha}^{(2n+1)}\right)\simeq SO_{2,1}(\RR)^{c\,(2n+1)}}:$ This group is connected and consists of all lifts from $D^{\eta}$ to $\hat M^{\eta}_{2n+1}$ of maps (\ref{groupetaalpha}).
\smallskip\\

Thus, we have proved the following theorem:

\begin{theorem}\label{main}\sl \hfill 

\noindent (i) A non-trivial cover of a hypersurface $\nu_{\alpha}$ with $-1<\alpha<1$  is $CR$-equivalent to either $\nu_{\alpha}^{(\infty)}$ or $\nu_{\alpha}^{(n)}$ with $n\ge 2$ defined in (\ref{universalnu}), (\ref{finitenu}); the groups of $CR$-automorphisms of these covers are described in {\bf (A)}-{\bf (B)}, in particular, for every $N\in\{2,3,\dots,\infty\}$ there exists a complex 2-dimensional manifold $M^{\nu}_N$ on which a connected 3-dimensional Lie group $G^{\nu}_N$ acts by holomorphic transformations, such that for every $-1<\alpha<1$ the hypersurface $\nu_{\alpha}^{(N)}$  is a $G^{\nu}_N$-orbit in $M^{\nu}_N$ and the group $\hbox{Aut}_{CR}\left(\nu_{\alpha}^{(N)}\right)$ consists of the restrictions of the elements of $G^{\nu}_N$ to $\nu_{\alpha}^{(N)}$.
\vspace{0.5cm}

\noindent (ii) A non-trivial cover of a hypersurface $\eta_{\alpha}$ with $\alpha>1$ is $CR$-equivalent to either $\eta_{\alpha}^{(\infty)}$ or $\eta_{\alpha}^{(n)}$ with $n\ge 2$ defined in (\ref{eta2}), (\ref{universaleta}), (\ref{finiteeta1}), (\ref{etanodd}); the groups of $CR$-automorphisms of the above covers are described in
{\bf (C)}-{\bf (E)}, in particular, for every $N\in\{2,3,\dots,\infty\}$ there exists a complex 2-dimensional manifold $M^{\eta}_N$ on which a connected 3-dimensional Lie group $G^{\eta}_N$ acts by holomorphic transformations, such that for every $\alpha>1$ the hypersurface $\eta_{\alpha}^{(N)}$  is a $G^{\eta}_N$-orbit in $M^{\eta}_N$ and the group $\hbox{Aut}_{CR}\left(\eta_{\alpha}^{(N)}\right)$ consists of the restrictions of the elements of $G^{\eta}_N$ to $\eta_{\alpha}^{(N)}$. 
\end{theorem}

{\obeylines
Department of Mathematics
The Australian National University
Canberra, ACT 0200
AUSTRALIA
E-mail: alexander.isaev@maths.anu.edu.au
}

\end{document}